\documentclass[12pt,intlim,righttag]{article}
\usepackage{amsmath,amsthm}
\usepackage{amssymb}
\usepackage{amsfonts}

\newcommand\beq{\begin{equation}}
\newcommand\eeq{\end{equation}}

\newcommand{\beaa}{\begin{eqnarray*}}
\newcommand{\eeaa}{\end{eqnarray*}}

\theoremstyle{Theorem}

\theoremstyle{corollary}

\theoremstyle{remark}

\theoremstyle{definition}

\begin{document}
\title{Functional equations and martingales
 }

\author{M. Mania$^{1)}$ and L. Tikanadze$^{2)}$}

\date{~}
\maketitle

\begin{center}
$^{1)}$ A. Razmadze Mathematical Institute of Tbilisi State University and
Georgian-American University,  Tbilisi, Georgia,
\newline(e-mail: misha.mania@gmail.com)
\\
$^{2)}$ Ivane Javakhishvili Tbilisi State University,  Tbilisi,
Georgia
\newline(e-mail: tikanadzeluka@gmail.com)
\end{center}

\begin{abstract}
{\bf Abstract.}We consider  functional equations (Cauchy's, Abel's  and some other functional equations) and show that to find general solution  of these equations  is equivalent to establish
that a space-transformation of a  Brownian Motion  by  suitable  function (or functions)  is a martingale.

\end{abstract}

\noindent {\it 2010 Mathematics Subject Classification. 60G44, 60J65, 97I70}

\noindent {\it Keywords}:  Martingales,  Functional Equations,  Brownian Motion

\section{Introduction}

The aim of this paper is to give a martingale characterization  of  the general measurable
solutions of  Cauchy's, Abel's and some other functional equations. We show that finding the general solution  of these equations  is equivalent to 
establishing that a space-transformation of a  Brownian Motion  by  a suitable  function (or functions)   is a martingale with almost surely right-continuous paths. 

A functional equation is an equation in which a function  (or a set of functions)  satisfying a certain relationship has to be found. The solution of functional equations is one of the oldest topics of mathematical analysis.
Such  equations  has 
applications  in many fields of  pure  mathematics as well, as  in applied science, such as   geometry, real and complex analysis,  partial differential equations,  probability theory, functional
analysis, dynamical systems,  decision analysis, economics, engineering and more. 

Although, differential equations provide powerful methods for solving functional equations, the differentiability assumptions are not directly required for functions accuring functional equations and in many applications weaker assumptions are needed.
That is what Hilbert dedicated the second part of his fifth problem, where he stated :" ...  In how far are the assertions which we can make in the case of  differentiable functions true under proper modifications  without this assumption?"  Motivated 
by this suggestion many researchers have treated various functional equations  with mild (or without any)  regularity assumptions.

A   fundamental equation in the theory of functional
equations is the Cauchy additive  functional equation
\begin{equation}\label{1}
 f(x + y) = f(x) + f(y),\;\;\;\; \text{for  all}\;\;\; x, y \in R
\end{equation}
and the related  three equations:

\begin{equation}\label{2}
 f(x + y) = f(x) f(y),\;\;\;\; \text{for  all}\;\;\;\;\; x, y \in R
\end{equation}

\begin{equation}\label{3}
 f(x ) + f(y) = f(xy),\;\;\;\;\;\; \text{for  all}\;\;\;\;\; x, y \in R_+
\end{equation}

\begin{equation}\label{4}
 f(x y) = f(x)  f(y),\;\;\;\;\;\;\;\;\;\;\;\; \text{for  all}\;\;\;\; x, y \in R_+
\end{equation}

 Equations  (\ref{2}),  (\ref{3}) and  (\ref{4}) are called Caushy's exponential, logarithmic and power functional equations, since the general
solutions of these equations are functions $e^{cx}$, $c\ln x$ and $x^c$ respectively, for some constant $c\in R$. Equations (\ref{2}-{4}) may be reduced to  the Cauchy additive functional equation (\ref{1}), or solved similarly.

The Cauchy functional equation (\ref{1}) has been investigated by many authors, under various "regularity"  conditions    and  each of them implies ( in the  case of real functions $ f : R \to R$), 
that  $f(x) = cx$ for some $c\in  R$. For instance, Cauchy \cite{CA} assumed that $f$  is continuous, Frechet \cite{F}, Banach \cite{BA} and  Sierpinski \cite{S}   showed  that the measurabelity of  $f$ is  sufficient. The most general result in this direction
( Kestelman \cite{KE}, Ostrowski \cite{O}) when $cx$ is the only solution of (\ref{1}) is an assumption on $f$  to be bounded  from one side on a measurable set of positive measure. Note that if the function is  Lebesgue measurable,
then it is bounded from one side on a measurable set of positive measure. On the other hand, Hamel \cite{H} investigated  equation (\ref{1}) without any  conditions on $f$ (with a use of the axiom of choice).  He showed that there exist also nonlinear solutions of 
(\ref{1})  and he found all such solutions. 

The functional equation (1) has been generalized or modified in many other directions. See \cite{A},\cite{KUC} or the recent paper \cite{R} for more details and related references. In this paper we consider only  Lebesgue measurable real functions.

The present paper was motivated by a note of  S. Smirnov \cite{SM}, where an application of Bernstein's characterization of the normal distribution  is given to show that any measurable solution of the Cauchy functional equation (1) is locally integrable.
We use this idea to show the integrability of the transformed processes $f(W_t)$, where $W=(W_t, t\ge 0)$ is a  Brownian Motion.

Let $W=(W_t, t\ge 0)$ be a standard Brownian Motion defined on a  probability space  $(\Omega, \cal F, P)$ with
filtration  
$F=({{\cal F}}_t,t\ge0)$  satisfying  the usual conditions of right-continuity and completeness.
A function $f=(f(x), x\in R)$ is called a semimartingale function of  the process $X$  if the transformed  process $(f(X_t), t\ge0)$ is 
a semimartingale. It was shown by Wang (\cite {W})  that every semimartingale function  
of Brownian Motion is locally difference of two convex functions.  More generally,  in  \cite{CI}  it was proved that for a given Markov process $X$ the process  $f(X_t)$ is a semimartingale if and only if  it is locally difference of two 
excesive functions.  In \cite{Ch}, \cite{MT}  the  description
 of time-dependent semimartingale functions of  Brownian Motion and diffusion processes   in terms of generalized derivatives was given. All these results imply that if $f(W_t)$ is a  right-continuous martingale,
then $f$ is a linear function. We use this fact several times in the paper and  for convenience   give  a direct  proof of this assertion   in  Theorem A1 of the Appendix. 

Our goal is  to relate functional equations with  semimartingale and martingale functions of Brownian motion and to give probabilistic proofs of some assertions on functional equations. We consider also stochastic versions of
the Cauchy functional equation.

Denote by $\cal M$ the class of martingales adapted to the filtration $F$ with $P$- almost surely right-continuous paths.

We show (Theorem 2.1 and Theorem 2.2) that the function $f=(f(x), x\in R)$  is a measurable solution of functional equation  (\ref{1}) (resp. (\ref{2}),  (\ref{3}),  (\ref{4}) ) if and only if the process
$f(W_t)$  (resp. $\ln f(W_t)$, $f(e^{W_t})$, $\ln f(e^{W_t})$) is a martingale from $\cal M$, zero at time zero.

We consider also stochastic versions of   Cauchy's functional equation (\ref{1})
  $$
f\left(x+W_{t}\right)=f\left(x\right)+f\left(W_{t}\right)\;\;\; \text{for  all}\;\;\;x\in R\;\;\; and\;\;\;t\ge 0,
$$
 $$
f\left(x+W_{1}\right)=f\left(x\right)+f\left(W_{1}\right)\;\;\; \text{for  all}\;\;\;x\in R 
$$
and show (Theorem 2.3) that the general measurable solutions of these equations coincide with  the general solution of equation (\ref{1}) $f(x)=cx$.

In section 4  we consider the Abel functional equation
\begin{equation}\label{abel}
f(x+y)=h(x-y)+g(xy), \;\;\;  \text{for  all}\;\;\;x\in R , y\in R,
\end{equation}
where $f,h,g:R\to R$ are real functions. In 1823 Abel \cite{AB} gave  differential solutions of this equation. The general solution of equation (\ref{abel}) was given by Aczel \cite{A} and by Lajko \cite{L} without any regularity conditions, in terms
of the additive  function. See also \cite{CHU},  where the general solution of (\ref{abel})  for a large class of fields was determined.

We show that to find the  general solution of Abels's equation  (\ref{abel}) is equivalent to find  general solution of a problem formulated in terms of martingales. In particular, we prove that (Theorem 4.1)
the triple $(f,h, g)$ is a measurable solution of  equation (\ref{abel}) if and only if

$K(W_t,y)\in\cal M$ for any $y\in R$,

$K(x, W_t)\in\cal M$ for any $x\in R$ and

$K(0,y)=K(x,0)= const$,

where the function $K$ is defined by
\begin{equation}\label{K}
K(x,y)= f(x+y) - h(x-y).
\end{equation}

The proof of this result is based  on:

 Theorem 3.1, where it was shown  that if a function $G=(G(x), x\in R)$ is a measurable solution of  the conditional Cauchy functional equation
\begin{equation}\label{cond}
G(x^2-y^2)=G(x^2)-G(y^2), \;\;\; \text{for  all}\;\;\;x\in R , y\in R
\end{equation}
then  the process $G(x+\sigma W_t)$ is a martingale for every $x, \sigma\in R$

 and on Theorem A2 from Appendix, which implies that the function $K$  defined  by (\ref{K}) should have the form
$$
K(xy)= axy+d,
$$
for some constants $a,d\in R$.

Finally, we give a probabilistic proof to establish general solution of the quadratic functional equation
$$
f(x+y)+f(x-y)=2f(x)+2f(y).
$$

\section{Cauchy's functional equations}

Let $W=(W_t, t\ge 0)$ be a standard Brownian Motion defined on a  complete probability space  $(\Omega, \cal F, P)$. Let
  $F=({{\cal F}}_t,t\ge0)$ be a filtration satisfying  the usual conditions of right-continuity and completeness. Assume that $F$ is larger than the filtration generated
by the Brownian Motion ${\cal F}^W_t=\sigma (W_s, s\le t)$ and that $W_t-W_s$ is independent of ${\cal F}_s$ whenever $0\le s\le t$, i.e.,  $(W_t,  t\ge 0)$ is
also a Brownian motion with respect to the filtration ${\cal F}_t$. Suppose also that the filtration $F$ is large enough to support an additional Brownian Motion $B$ 
independent of $W$.

{\bf Theorem 2.1}
Let $f=\left(f\left(x\right),x\in R\right)$ be function of one variable. The following assertions are equivalent : 

$(a)$ The function  $f=\left(f\left(x\right),x\in R\right)$ is a measurable solution of the Cauchy functional equation
\ \begin{equation}\label{c1}
f(x+y)=f(x)+f(y)\;\;\;\; \text{for all}\;\;\;\;\; x,y\in R.
\end{equation}

$(b)$  $f=\left(f\left(x\right),x\in R\right)$ is a measurable function such that for any fixed  $t\ge0$
\begin{equation}\label{c2}
 f(x+W_{t})=f(x)+f(W_{t})
\end{equation} 
$P-$ a.s. for all $x\in R$. 

$(c)$  The process $\left(f\left(W_{t}\right),t\geq 0\right)\in\cal M$, i.e., it is a martingale  with  $P$- a.s. right-continuous paths, zero at time zero. 

 $(d)$   $f\left(x\right)=cx$ for some constant $c\in R$.

\begin{proof}
The proof of $(a)\to(b),  (d)\to(a)$ is evident and the proof of implication $(c)\to (d)$ follows from  Theorem A1 of Appendix. Let us show the implication $(b)\to(c)$.

 Let first show that
 $f(W_t)$ is ${\cal F}_t$-measurable for every $t\ge 0$. It is well known that there exists a Borel measurable function $\widetilde f$ such that $L(x: f(x)\neq \widetilde f(x))=0$, where
    $L$ is the Lebesgue measure. Then 
 $$
    P(\omega : f(W_t)\neq \widetilde f(W_t))=\int_R I_{(x: f(x)\neq \widetilde f(x))}\frac{1}{\sqrt{2\pi t}}e^{-\frac{x^2}{2t}}dx=0
  $$
    and $f(W_t)$ and ${\widetilde f}(W_t)$ are equivalent. Since ${\tilde f}(W_t)$ is ${\cal F}_t$- measurable and ${\cal F}_t$ is completed with $P$- null sets from ${\cal F}$,
    $f(W_t)$ will be also ${\cal F}_t$-measurable.

To show that $f\left(W_{t}\right)$ is integrable for any $t\geq 0$  we shall use the idea from \cite{SM} on application of the Bernstein theorem.

 Let $B_{t}$ be a standard Brownian motion independent of $W_{t}.$ and let
$$
 X=f\left(W_{t}\right),\;\;\;\;Y=f\left(B_{t}\right).
$$
 It follows from $(\ref{c2})$ that  $P$- a.s.
$$
X+Y=f\left(W_{t}\right)+f\left(B_{t}\right)=f\left(W_{t}+B_{t}\right)
$$
and substituting $x=B_{t}-W_{t}$ in  $(\ref{c2})$ we have 
$$
Y-X=f\left(B_{t}\right)-f\left(W_{t}\right)=f\left(B_{t}-W_{t}\right)\;\;\;\;P-a.s.
$$

Since $B_{t}+W_{t}$ and $B_{t}-W_{t}$ are independent, the random variables $f(B_{t}+W_{t})$ and $f(B_{t}- W_{t})$ will be also independent. 
Therefore, Bernstein's theorem (see Theorem A3 from Appendix) implies that $f(W_{t})$ (and $f(B_{t})$ also) is distributed normally. Hence 
$$
E|f(W_{t})|<\infty.
$$

Note that $f(W_{t})$ is integrable  also at any power.

Let us show now the martingale equality 
$$
E(f(W_{t})|\mathcal{F}_{s})=f(W_{s}),\;\;\;\;P- \text{a.s}
$$
for all $s,t ( s\leq t)$.

Substituting $x=W_{t}-W_{s}$ in the equality 
$$
f\left(x+W_{s}\right)=f\left(x\right)+f\left(W_{s}\right)
$$
we have that $P$-a.s.
\begin{equation}\label{c3}
f\left(W_{t}\right)-f\left(W_{s}\right)=f\left(W_{t}-W_{s}\right).
\end{equation}
Interchanging t and s in (\ref{c3})
\begin{equation}\label{c4}
f\left(W_{s}\right)-f\left(W_{t}\right)=f\left(W_{s}-W_{t}\right)
\end{equation}
and from  (\ref{c3}),  (\ref{c4}) we get
\begin{equation}
f\left(W_{t}-W_{s}\right)=-f\left(W_{s}-W_{t}\right)\;\;\;\;P-a.s.
\end{equation}
This implies that 
\begin{equation}\label{c5}
Ef\left(W_{t}-W_{s}\right)=0
\end{equation}
since $f\left(W_{t}-W_{s}\right)$ and $f\left(W_{s}-W_{t}\right)$ have the same distributions.

Taking conditional expectations in  (\ref{c3}), since $f\left(W_{t}-W_{s}\right)$ is independent of $\mathcal{F}_{s},$ we obtain 
$$
E\left(f\left(W_{t}\right)-f\left(W_{s}\right)|\mathcal{F}_{s}\right)=E\left(f\left(W_{t}-W_{s}\right)|\mathcal{F}_{s}\right)=Ef\left(W_{t}-W_{s}\right)\;\;\;\;\;P-\text{a.s.}
$$
Therefore (\ref{c5}) implies that $P$-a.s
$$
E\left(f\left(W_{t}\right)-f\left(W_{s}\right)|\mathcal{F}_{s}\right)=0,
$$
hence $\left(f\left(W_{t}\right),\mathcal{F}_{t},t\geq 0\right)$ is a martingale.

It follows from equality (\ref{c2})  that $f(x+W_t)$ is also a martingale for every $x\in R$, which implies 
(implication $c)\to a)$ of Theorem A1) that $f(W_t)$  is a  $P$- a.s.  right-continuous  martingale.

{\bf Remark.} Note that, if almost all paths of the process $f(W_t)$ are  right-continuous, then the function $f(x)$ is continuous. Thus, the continuity and 
right-continuity for the transformed process $f(W_t)$ are equivalent.
\end{proof} 

\

Now let us consider Cauchy's remaining three functional equations.  Denote by $R_+$ the set of positive numbers.

{\bf Theorem 2.2}

$(a)$  The function $\left(f\left(x\right),x\in R\right)$ is a measurable non-zero solution of functional equation 
\begin{equation}\label{ex}
f\left(x+y\right)=f\left(x\right)f\left(y\right), \;\;\;\; x,y\in R
\end{equation}
 if and only if $f\left(W_{t}\right)$ is strictly positive process such that $\ln f\left(W_{t}\right)\in\cal M$, i.e., is a martingale with  $P$- a.s. right-continuous paths.
 zero at time zero.
.

(b)   The function $\left(f\left(x\right),x\in R\right)$ is a measurable solution of functional equation 
\begin{equation}\label{lg}
f\left(x\right)+f\left(y\right)=f\left(xy\right), \;\;\;\;x,y\in R_{+} 
\end{equation}

if and only if the process $f\left(e^{W_{t}}\right)\in\cal M$ and equals to zero at time $t=0$. 
.

 (c) The function $\left(f\left(x\right),x \in R\right)$ is a measurable non-zero solution of functional equation 
\begin{equation}\label{pw}
f\left(xy\right)=f\left(x\right)f\left(y\right), \;\;\;\;x,y \in R_{+}
\end{equation}
if and only if $f\left(e^{W_{t}}\right)$ is a strictly positive process such that  $\ln f\left(e^{W_{t}}\right)\in\cal M$ and equals to zero at time $t=0$.

\begin{proof}
We shall prove assertion (c). The proofs of (a) and  (b) are similar. It is obvious (and well known) that a solution of (\ref{pw}) is either everywhere or nowhere 0. 

Indeed, (\ref{pw}) implies that \[
\begin{split}
f\left(x^{2}\right)=f^{2}\left(x\right)\geq 0
\end{split}
\]
and if $f\left(x_{0}\right)=0$ for some $x_{0}>0$ then \[
\begin{split}
f\left(x\right)=f\left(x_{0}\frac{x}{x_{0}}\right)=f\left(x_{0}\right)f\left(\frac{x}{x_{0}}\right)=0.
\end{split}
\]
Therefore, excluding the solution $f\left(x\right)=0$ for all $x>0$ we will have that $f\left(x\right)>0$ for all $x>0$ and the process $f\left(e^{W_{t}}\right)$ will be strictly positive. 

Let us show that the process $\left(lnf\left(e^{W_{t}}\right),t\geq 0\right)$ is a martingale. Let first show that \[
\begin{split}
E|lnf\left(e^{W_{t}}\right)|<\infty
\end{split}
\]
for all $t\geq 0$.\\Let $X=f\left(e^{W_{t}}\right)$ and $Y=f\left(e^{B_{t}}\right)$, where $B_{t}$ is a Brownian motion independent of $W_{t}$. It follows from (\ref{pw}) that 
\begin{equation}\label{xy}
XY=f\left(e^{W_{t}}\right)f\left(e^{B_{t}}\right)=f\left(e^{W_{t}+B_{t}}\right),
\end{equation}
\begin{equation}\label{xy1}
 \frac{X}{Y}=\frac{f\left(e^{W_{t}}\right)}{f\left(e^{B_{t}}\right)}=f\left(e^{W_{t}-B_{t}}\right)
\end{equation}

Since $W_{t}+B_{t}$ and $W_{t}-B_{t}$ are independent,  it follows from equations (\ref{xy}) and (\ref{xy1}) that the random variables $XY$ and $\frac{X}{Y}$ will be also independent. 
Therefore, by Bernstein's theorem $X=f\left(e^{W_{t}}\right)$ (and $Y=f\left(e^{B_{t}}\right)$) 
will have the log-normal distribution and $lnf\left(e^{W_{t}}\right)$ admits the normal distribution, hence $lnf\left(e^{W_{t}}\right)$ is integrable for any $t\geq 0$. 

By change of variables and functions the equation (\ref{pw}) goes over into  \[
\begin{split}
f\left(e^{u+v}\right)=f\left(e^{u}\right)f\left(e^{v}\right)
\end{split}
\] 
and substituting $u=W_{t}-W_{s}$ and $v=W_{s}$ in this equation and taking logarithms we have that 
\begin{equation}\label{ln}
lnf\left(e^{W_{t}}\right)-lnf\left(e^{W_{s}}\right)=lnf\left(e^{W_{t}-W_{s}}\right)
\end{equation}
By independent increment of Brownian motion $lnf\left(e^{W_{t}-W_{s}}\right)$ is independent of $\mathcal{F}_{s}$ and taking conditional expectation in  (\ref{ln}) we have that
$P-a.s.$ 
\[
\begin{split}
E\left(lnf\left(e^{W_{t}}\right)-lnf\left(e^{W_{s}}\right)|\mathcal{F}_{s}\right)=E\left(lnf\left(e^{W_{t}-W_{s}}\right)|\mathcal{F}_{s}\right)=Elnf\left(e^{W_{t}-W_{s}}\right)
\end{split}
\]
But \[
\begin{split}
Elnf\left(e^{W_{t}-W_{s}}\right)=0,
\end{split}
\]
since the function $lnf\left(e^{u}\right)$ is odd and the distribution of $W_{t}-W_{s}$ is symmetric.

Since
$$
\ln f(e^{x+W_t})= \ln f(e^{x})+\ln f(e^{W_t}),
$$
the process $\ln f(e^{x+W_t})$ will be also a martingale for any $x\in R$, which implies (see Theorem A1 implication $c)\to a)$) that  almost all paths of the martingale $\ln f(e^{W_t})$ are  right-continuous.

Now let us assume that process $lnf\left(e^{W_{t}}\right)$ is a  $P$-a.s. right-continuous martingale, zero at time zero. Then Theorem A1 (implication $a)\to b)$) implies that \[
\begin{split}
lnf\left(e^{u}\right)=\lambda u
\end{split}
\]
for some $\lambda \in R$ and changing variables $u=lny$ we obtain that $f\left(y\right)=y^{\lambda}$, which satisfies equation (\ref{pw}).
\end{proof}

Now let us show that if equality (\ref{c2}) is satisfied only for  $t=1$ the set of solutions remains as it was. I.e., we consider the following stochastic version of Cauchy's functional equation 
  \begin{equation}\label{xi}
f\left(x+\xi\right)= f\left(x\right)+f\left(\xi\right)\;\;\;\;\;\;\;\;\;\;P\text{- a.s.}\;\;\text{ for all}\;\;\; x\in R,
  \end{equation}
 where $\xi$ is a random variable with standard normal distribution , i.e. \[
  \begin{split}
 E\xi=0, \quad E\xi^{2}=1
 \end{split}
 \]
 The following theorem shows that  (\ref{xi}) is also equivalent to assertions $(a)- (d)$  of Theorem 2.1

{\bf Theorem 2.3}
Any measurable solution of (\ref{xi}) is linear.

\begin{proof}
It is evident that if  $f$  is a solution of (\ref{xi}), then $f\left(0\right)=0$ and substituting $x=-\xi$ in (\ref{xi}) we have that $f\left(\xi\right)=-f\left(-\xi\right)$  $P-a.s.$, which implies that \[
\begin{split}
Ef\left(\xi\right)=0.
\end{split}
\]
since $\xi$ is symmetrically distributed.
Similarly as in Theorem 2.1, one can show that the random variable $f\left(\xi\right)$ is also normally distributed. This implies that $f\left(\xi\right)$ is square integrable \[
\begin{split}
Ef^{2}\left(\xi\right)=\frac{1}{\sqrt{2\pi}}\int_{R}{f^{2}\left(x\right)e^{-\frac{x^{2}}{2}}dx}<\infty
\end{split}
\]
and the function $f\left(x\right)$ is locally square integrable.Taking expectation in (\ref{xi}) we obtain that \[
\begin{split}
f\left(x\right)=Ef\left(x+\xi\right)=\int_{R}{f\left(x+y\right)\frac{1}{\sqrt{2\pi}}e^{-\frac{y^{2}}{2}}dy} 
\end{split}
\]
and after changing variables $x+y=z$ we get \[
\begin{split}
f\left(x\right)=\int_{R}{f\left(z\right)\frac{1}{\sqrt{2\pi}}e^{-\frac{\left(z-x\right)^{2}}{2}}dz}.
\end{split}
\]
It follows from here  that $f\left(x\right)$ is differentiable and 
$$
f^{\prime}\left(x\right)=\int_{R}{f\left(z\right)\left(z-x\right)\frac{1}{\sqrt{2\pi}}e^{-\frac{\left(z-x\right)^{2}}{2}}dz}=
$$
\begin{equation}\label{xi2}
=\int_{R}{f\left(x+y\right)y\frac{1}{\sqrt{2\pi}}e^{-\frac{y^{2}}{2}}dy}=Ef\left(x+\xi\right)\xi .
\end{equation}
Using (\ref{xi}), (\ref{xi2})  and equality $E\xi=0$ , we obtain that 
\begin{equation}\label{xi3}
f^{\prime}\left(x\right)=Ef\left(x+\xi\right)\xi=Ef\left(x\right)\xi+E\xi f\left(\xi\right)=E\xi f\left(\xi\right)
\end{equation}
Note that $\xi f\left(\xi\right)$ is integrable, since $\xi$ and $f\left(\xi\right)$ are Gaussian and hence square integrable. \\Thus, (\ref{xi3}) implies that $f^{\prime}\left(x\right)$ is constant and $f\left(x\right)=\lambda x$ for some $\lambda \in R$.
\end{proof}

\section {Cauchy conditional functional equation}

Let consider the conditional Cauchy functional  equation
\begin{equation}\label{cond2}
G(x^2-y^2)=G(x^2)-G(y^2), \;\;\; \text{for  all}\;\;\;x\in R , y\in R.
\end{equation}
It is well known (see e.g. \cite{CHU}) that $G$ is an additive map.
We give an equivalent formulation in terms of corresponding martingale problem.

Let first mention some simple properties of equation (\ref{cond2}) which will be used in the sequel.
It is evident that
$$
G(0)=0\;\;\;\;\text{and}\;\;\;\;G(u)=- G(-u).
$$
Since for any $x,y\in R$ there exists $z\in R$ such that $x^2+y^2=z^2$, it follows from (\ref{cond2}) that
$$
G(x^2)=G(z^2-y^2)=G(z^2)-G(y^2)
$$
and hence
\begin{equation}\label{g2}
G(x^2+y^2)=G(x^2)+G(y^2), \;\;\;  \text{for  all}\;\;\;x\in R ,\;\;y\in R.
\end{equation}

{\bf Theorem 3.1}. The function $G=(G(x), x\in R)$ is a measurable solution of (\ref{cond2}) if and only if  the process $(G(W_t), t\ge 0)$ is a  $P$- a.s.  right-continuous  martingale, zero at time zero.

{\it Proof}. 
Assume that  $G=(G(x), x\in R)$ is a measurable solution of (\ref{cond2}). Let us show that 
$G(x+\sigma W_t)$ is a martingale for any $x, \sigma\in R$.

Let $\xi^+=\max(\xi, 0)$ and $\xi^-=-\min(\xi, 0)$ be the positive and negative parts of random variable $\xi$. 

Let 
\begin{equation}\label{not}
M_t\equiv x+\sigma W_t.
\end{equation}

 Since

\begin{equation}\label{g3}
M_t= M_t^+ - M_t^-,
\end{equation}
it follows from (\ref{cond2}) that
$$
G(M_t)= G(M_t^+ -  M_t^-)= G(\big(\sqrt{M_t^+}\big)^2 - (\sqrt{M_t^-})^2)=
$$
\begin{equation}\label{g4}
 G(\big(\sqrt{M_t^+}\big)^2) - G(\big (\sqrt{M_t^-}\big)^2)= G(M_t^+) - G(  M_t^-).
\end{equation}

Therefore, $P-a.s$ 
$$
E(G(M_t)-G(M_s)/F_s)=
$$
$$
=E(G(M_t^+)-G(M_t^-)- G(M_s^+)+G(M_s^-)/F_s)= \;\;\;\;\;\;\;(\text{by}\;\;(\ref{g4}))
$$
$$
=E(G(M_t^+ +M_s^-)- G(M_t^- +M_s^+)/F_s)= \;\;\;\;\;\;\;(\text{by}\;\;(\ref{g2}))
$$
$$
=E(G(M_t^+ +M_s^- -M_t^- -M_s^+)/F_s)= \;\;\;\;\;\;\;(\text{by}\;\;(\ref{cond2}))
$$
$$
=E(G(M_t - M_s)/F_s)=\;\;\;\;\;\;\;(\text{by}\;\;(\ref{g3}))
$$
$$
=EG(M_t-M_s)\;\;\;\;\;(\text{by independent increments of the  BM)}
$$
$$
=E G(\sigma W_t- \sigma W_s)\;\;\;\; (\text{by notation}  \;\; (\ref{not})).
$$
Finally, $EG(\sigma W_t- \sigma W_s)=0$, since $(\sigma W_t- \sigma W_s)^+$ and $(\sigma W_t- \sigma W_s)^-$ are identically distributed and by (\ref{g4})
$$
EG(\sigma W_t- \sigma W_s)=EG(\sigma W_t- \sigma W_s)^+- EG(\sigma W_t- \sigma W_s)^-=0.
$$
To show the integrability of $G(M_t)$ we shall use again the Bernstein theorem. Let $B_t$ be a Brownian Motion independent of $W_t$ and let
$$
X=G(x+\sigma W_t)\;\;\;\;\text{and}\;\;\;Y=G(x+\sigma B_t).
$$
Then using successively equations  (\ref{g4}), (\ref{g2}), (\ref{cond2}) and (\ref{g3}) we obtain that
$$
X-Y=G(x+\sigma W_t)-G(x+\sigma B_t)=G((\sigma (W_t- B_t)),
$$
$$
X+Y=G(x+\sigma W_t)+G(x+\sigma B_t)=G(\sigma (W_t+ B_t)).
$$
Since $W_t+B_t$ and $W_t-B_t$ are independent, the random variables  $G(\sigma (W_t+ B_t))$ and    $G((\sigma (W_t- B_t))$ will be also independent and by Bernstein's theorem 
$G(x+\sigma W_t)$ will have the normal distribution, which implies that  $E|G(x+\sigma W_t)|<\infty$ for every $t\ge 0$ and for all $x,\sigma\in R$. Thus, the process $G(x+\sigma W_t)$ is a 
martingale, which implies (implication $c)\to a)$ of Theorem A1) that $G(\sigma W_t)$  is a  $P$- a.s.  right-continuous  martingale, zero at time zero.

Vice-versa, if $G(W_t)\in\cal M$ and $G(0)=0$, then by Theorem A1 from Appendix  $G(x)=\lambda x$ for some constant $\lambda\in R$  and
it is evident that $G(x)=\lambda x$ satisfies equation (\ref{cond2}).

\section {The Abel  functional equation}
In 1823 Abel \cite{AB} considered functional equation
\begin{equation}\label{abel2}
f(x+y)=h(x-y)+g(xy), \;\;\;  \text{for  all}\;\;\;x\in R , y\in R,
\end{equation}
where $f,h,g:R\to R$ are real functions. In the same manuscript  Abel \cite{AB} gave  differential solutions of this equation. The general solution 
of equation (\ref{abel2}) was given by Aczel \cite{A} and by Lajko \cite{L} without any regularity conditions, in terms  of additive map.

We show that to find the  general solution of Abel's equation  (\ref{abel2}) is equivalent to find  the general solution of a problem formulated in terms of martingales. 

Let us define the function
$$
K(x,y)= f(x+y) - h(x-y).
$$

{\bf Theorem 4.1}
The triple $(f,h, g)$ is a measurable solution of  Abel's  functional  equation (\ref{abel2}) if and only if

$K(W_t,y)\in\cal M$  for any $y\in R$,

$K(x, W_t)\in\cal M$ for any $x\in R$ 

and  $K(0,y)=K(x,0)= \lambda,$

where  $ \lambda$ is some constant.

{\it{Proof.}}
Let the triple $(f, h, g)$ be a measurable solution of (\ref{abel2}). Let
$$
H(x)=h(x) -h(0)\;\;\;\text{and}\;\;\;\;G(x)= g(x) - g(0).
$$
Then it is easy to see that the pair $(H,G)$ satisfies the functional equation
\begin{equation}\label{HG}
H(x+y)-H(x-y)=G(xy), \;\;\;  \text{for  all}\;\;\;x\in R , y\in R.
\end{equation}
Indeed, from (\ref{abel2}) taking $y=0$ we have that $f(x)=h(x)+g(0)$. Therefore
\begin{equation}\label{a3}
f(x+y)=h(x+y)+g(0)=H(x+y)+h(0)+g(0),
\end{equation}
\begin{equation}\label{a4}
h(x-y)= H(x-y)+h(0)
\end{equation}
and  (\ref{abel2}),  (\ref{a3}) and  (\ref{a4}) imply that
\begin{equation}\label{a5}
G(xy)=g(xy)-g(0)=f(x+y)-h(x-y)-g(0)=H(x+y)- H(x-y),
\end{equation}
hence the pair $(H, G)$ satisfies (\ref{HG}).

It follows from (\ref{HG}) that the function $G$ satisfies the conditional Cauchy equation (\ref{cond2}), since
\begin{equation}\label{a6}
G(u^2-v^2)= G((u+v)(u-v))=H(2u)- H(2v)= G(u^2)- G(v^2)
\end{equation}
for all $u,v\in R$.
According to the proof of Theorem 3.1 the process $G(\sigma W_t)\in\cal M$  for any $\sigma\in R$. Since by (\ref{abel2}) 
$$
K(x,y)= f(x+y) - h(x-y)=g(xy)= G(xy)+g(0),
$$
we obtain that
$$
K(W_t, y)=G(yW_t)+g(0)\in\cal M\;\;\text{and}
$$
$$
K(x, W_t)=G(xW_t)+g(0)\in\cal M.
$$
It is evident that 
$$
K(0,y)=K(x,0)=g(0)\equiv\lambda.
$$

Now let us assume that 
$K(W_t,y)$ and $K(x, W_t)$ belong to $\cal M$  for any $y\in R$ and $x\in R$ respectively, with
$K(0,y)=K(x,0)= \lambda.$

It follows from  Theorem A2 of the Appendix, that $K(x,y)$ will be of the form
$$
K(x,y)=axy+bx+cy+d.
$$
Condition $K(0,y)=K(x,0)= \lambda$  implies that
$$
bx+d=\lambda \;\;\;\;\text{and}\;\;\;\;cy+d=\lambda
$$
for all $x,y\in R$. Hence $b=c=0$ and
$$
K(x,y)=axy+d.
$$
Thus,
\begin{equation}\label{axy}
f(x+y)-h(x-y)=axy+d.
\end{equation}
Taking $x=y=\frac{u}{2}$  in(\ref{axy}) we have

\begin{equation}\label{axy2}
f(u)=\frac{a}{4}u^2+h(0)+d
\end{equation}
and if we take $x=u, y=0$ we obtain from (\ref{axy}) and  (\ref{axy2}) that
\begin{equation}\label{axy3}
h(u)=f(u)-d=\frac{a}{4}u^2+h(0).
\end{equation}
Therefore, it follows from  (\ref{axy2}) and  (\ref{axy3}) that
$$
f(u+v)-h(u-v)=
$$
$$
=\frac{a}{4}(u+v)^2+h(0)+d-\frac{a}{4}(u-v)^2-h(0)
$$
$$
=auv+d.
$$
Hence, the triple
$$
g(x)=ax+d,
$$
$$
h(x)=\frac{a}{4}x^2+h(0),
$$
$$
f(x)=\frac{a}{4}x^2+h(0)+d,
$$
where  $a, d$ and $ h(0)$ are constans, satisfies equation (\ref{abel2}). This proves also that it gives the general solution of (\ref{abel2}).

\section{Quadratic functional equations}

Let us consider quadratic functional equation 
\begin{equation}\label{q1}
f\left(x+y\right)+f\left(x-y\right)=2f\left(x\right)+2f\left(y\right)
\end{equation}
for all $x,y \in R$. It is well known (see, e.g., \cite{KA}), that the general solution of equation $(\ref{q1}))$ is the function $f\left(x\right)=\lambda x^{2}$. In the following theorem we give a probabilistic proof of this assertion.

{\bf Theorem 5.1}
The general measurable solution of equation $(\ref{q1})$ is of the form 
\begin{equation}\label{q2}
f\left(x\right)=\lambda x^{2}
\end{equation}
where $\lambda\in R$ is some constant.

\begin{proof}
It is evident that if $f$ is a solution of  $(\ref{q1})$ then $f\left(0\right)=0$ and \[
\begin{split}
f\left(x\right)=f\left(-x\right), \quad \text{for all} \quad x\in R.
\end{split}
\]
Let 
\begin{equation}\label{q3}
G\left(x,y\right)=f\left(x+y\right)-f\left(x\right)-f\left(y\right).
\end{equation}
It is easy yo see that 
\begin{equation}\label{q0}
G\left(0,x\right)=G\left(y,0\right)=0 \quad \text{and}
\end{equation}

\begin{equation}\label{q4}
G\left(x,y\right)=-G\left(-x,y\right)=-G\left(x,-y\right)
\end{equation}
Let us show that the process ($G\left(x,W_{t}\right),t\geq 0$) is a martingale for any $x\in R$ and ($G\left(W_{t},y\right),t \geq 0$) is a martingale for any $y\in R$. 

After simple transformations it follows from  $(\ref{q1})$,   $(\ref{q3})$ and  the equality  $f\left(W_{t}-2W_{s}+y\right)=f\left(2W_{s}-W_{t}-y\right)$ that 

\[
\begin{split}
G\left(W_{t},y\right)=f\left(W_{t}+y\right)-f\left(W_{t}\right)-f\left(y\right)=\\ \frac{1}{2}\left[f\left(W_{t}+y\right)-f\left(W_{t}-y\right)\right]=\\ \frac{1}{2}\left[f\left(W_{t}-W_{s}+y+W_{s}\right)-f\left(W_{s}-y+W_{t}-W_{s}\right)\right]=\\ \frac{1}{2}[2f\left(W_{t}-W_{s}+y\right)+2f\left(W_{s}\right)-f\left(W_{t}-2W_{s}-W_{s}\right)-\\-2f\left(W_{s}-y\right)-2f\left(W_{t}-W_{s}\right)+f\left(2W_{s}-W_{t}-y\right)]=\\ f\left(W_{t}-W_{s}+y\right)-f\left(W_{t}-W_{s}\right)-f\left(y\right)+\\+ f\left(W_{s}\right)+f\left(y\right)-f\left(W_{s}-y\right)=\\ G\left(W_{t}-W_{s},y\right)+f\left(W_{s}+y\right)-f\left(W_{s}\right)-f\left(y\right)=\\ =G\left(W_{t}-W_{s},y\right)+G\left(W_{s},y\right).
\end{split}
\]
Thus,
\begin{equation}\label{q5}
G\left(W_{t},y\right)-G\left(W_{s},y\right)=G\left(W_{t}-W_{s},y\right)
\end{equation}
and taking conditional expectations  in (\ref{q5}) we get  that $P$-a.s.
$$
E\left(G\left(W_{t},y\right)-G\left(W_{s},y\right)|\mathcal{F}_{s}\right)=
$$ 
$$
= E\left(G\left(W_{t}-W_{s},y\right)|\mathcal{F}_{s}\right)=EG\left(W_{t}-W_{s},y\right)=0.
$$
Here we used the independent increment property of the Brownian motion (hence $G(W_{t}-W_{s},y)$ is independent of $\mathcal{F}_{s}$), the symmetric distribution  of $W_{t}-W_{s}$ and 
that $G(x, y)$ is odd for any  $y$.

Similarly one can show that the processes $G(a+W_t, y)$ and $G(x, b+W_{t})$ are martingale for any $a,b\in R$ respectively and by Theorem A1  $G(W_t, y)\in\cal M$ and $G(x, W_t)\in\cal M$  for any $x, y\in R$.

Therefore, it follows from Theorem A2 of the   Appendix (taking (\ref{q0}) in mind)
 that 
$$
G\left(x,y\right)=axy
$$
for some constant $a\in R$. Finally,  from $(\ref{q3}), $ taking $y=x$ we obtain that 
$$
ax^{2}=G\left(x,x\right)=f\left(2x\right)-2f\left(x\right)=2f\left(x\right),
$$

hence $f\left(x\right)=\frac{a}{2}x^{2}$.
\end{proof}

\section{Appendix}

{\bf Theorem A1.}  Let  $\left(f\left(x\right),x\in R\right)$ be a  function of one variable.  The following assertions are equivalent.:

$a)$  $f(W_t), t\ge0\in\cal M$, i.e., it is a martingale with  $P$- a.s.  right-continuous paths.

$b)$ The function  $f$ is linear,
\begin{equation}\label{abx}
f\left(x\right)=ax+b
\end{equation}
for some constants $a, b\in R$. 

$c)$  The process $f(x+W_t), t\ge0$ is a martingale for every $x\in R$.
\begin{proof}
$a)\to b)$
Let $\left(f\left(W_{t}\right),\mathcal{F}_{t},t\geq 0\right)$ be a martingale with  $P$- a.s.  right-continuous paths. Then the function $\left(f\left(x\right),x\in R\right)$
will be continuous. Let \[
\begin{split}
g\left(t,x\right)=E\left(f\left(W_{T}\right)|W_{t}=x\right).
\end{split}
\]
It is well known that $g\left(t,x\right)$ satisfies the Backward Kolmogorov's equation \[
\begin{split}
\frac{\partial g}{\partial t}+\frac{1}{2}\frac{\partial^{2}g}{\partial x^{2}}=0.
\end{split}
\]
By the Markov property of the Brownian motion \[
\begin{split}
g\left(t,W_{t}\right)=E\left(f\left(W_{T}\right)|\mathcal{F}_{t}\right)\;\;\;\;\text{a.s.}
\end{split}
\]
and from the martingale property of $f\left(W_{t}\right)$ we have that for all $t\le T$ \[
\begin{split}
g\left(t,W_{t}\right)=f\left(W_{t}\right) \quad  a.s.
\end{split}
\]
Therefore, for all $t\le T$ \[
\begin{split}
\int_{}{|g\left(t,x\right)-f\left(x\right)|\frac{1}{\sqrt{2\pi t}}e^{-\frac{x^{2}}{2t}}dx}=0
\end{split}
\]
which implies that for any $t\le T$ \[
\begin{split}
g\left(t,x\right)=f\left(x\right) \quad a.e
\end{split}
\]
with respect to the Lebesgue measure. Since $T$ is arbitrary, by continuity of $f$ and  $g$ \[
\begin{split}
g\left(t,x\right)=f\left(x\right)
\end{split}
\]
for any $t>0$.\\ Thus $g\left(t,x\right)$ does not depend on t and $\frac{\partial g}{\partial t}=0.$ Therefore \[
\begin{split}
\frac{\partial^{2}g\left(t,x\right)}{\partial x^{2}}=\frac{\partial ^{2}f}{\partial x^{2}}=0,
\end{split}
\]
which implies that $f\left(x\right)$ is of the form (\ref{abx}).

The implication $b)\to c)$  is evident.

To prove the implication $c)\to a)$ we  note that by  the martingale equality  we have that
$$
f(x)=E f(x+W_t)=\int_R f(x+y) \frac{1}{\sqrt{2\pi t}}e^{-\frac{y^2}{2t}}dy=
$$
 \begin{equation}\label{cont}
=\int_R f(y) \frac{1}{\sqrt{2\pi t}}e^{-\frac{(y-x)^2}{2t}}dy.
\end{equation}
Since $E|f(x+W_t)|<\infty$, equality (\ref{cont}) implies that the function  $f(x)$ is continuous. Because almost all paths of Brownian Motion are continuous,  the process  $g(W_t)$
will be continuous  $P$- a.s.
\end{proof} 

{\bf Remark.}  If the transformed process $f(W_t)$ is a ${\cal F}_t$-martingale, then it will be a martingale with respect to the natural filtration ${\cal F}^W$. 
Therefore $f(W_t)$, as any ${\cal F}^W$-martingale, will have a continuous modification, but itself it can be not continuous. If $f(W_t)$ is only a martingale 
 (without assuming the regularity of paths), then $f(x)$  will coincide with a linear function almost everywhere with respect to the Lebesgue measure.

{\bf Theorem A2.} A function $G=(G(x, y), x, y\in R)$ is of the form
\begin{equation}\label{form}
G(x,y)=axy+bx+cy+d,
\end{equation}
where $a, b, c $ and $d$ are some constants,if and only if

$G(W_t,y)\in\cal M$ for any $y\in R$ and

$G(x, W_t)\in \cal M$  for any $x\in R$.

\begin{proof}
If $G(W_t,y)\in\cal M$ for any $y\in R$, it follows from Theorem A1 that
\begin{equation}\label{form2}
G(x,y)=\alpha(y) x+\beta(y).
\end{equation}
Since  $G(x, W_t)\in\cal M$ for any $x\in R$, the process
\begin{equation}\label{form3}
\alpha(W_t) x+\beta(W_t)
\end{equation}
will be a martingale from $\cal M$ for any $x\in R$, which implies that 
$\alpha(W_t)\in\cal M$, $\beta(W_t)\in\cal M$  and using again Theorem A1 we have that
\begin{equation}\label{form4}
\alpha(y)= ay+b\;\;\;\;\text{and}\;\;\;\;\beta(y)= cy+d
\end{equation}
for some constants $a, b, c$ and $d$.

Therefore, substituting expressions of $\alpha(y)$ and $\beta(y)$ in (\ref{form2}) we obtain the representation (\ref{form}).

The inverse assertion is obvious.
\end{proof}

\

The following result was proved by Bernstein \cite{BE} under assumption of equal (and finite) variances. We shall use general version of Bernstein's theorem due to Quine \cite{Q}.

{\bf Theorem A3.}
Assume that X and Y are independent random variables. Let $Z=X+Y$ and $V=X-Y$. If Z and V are independent, then X and Y are normally distributed with the same variances.

\end{document}